\documentclass[12pt]{article}
\usepackage{latexsym}
\usepackage{amsmath,epsfig}
\usepackage{amssymb}
\usepackage{amsfonts}

\numberwithin{equation}{section}

\title{Recovery of time dependent volatility coefficient by
 linearization}

\author{ Victor Isakov}

\date{July 2013}

\begin{document}
\maketitle
\newtheorem{theorem}{Theorem}[section]
\newtheorem{lemma}[theorem]{Lemma}
\newtheorem{corollary}[theorem]{Corollary}
\renewcommand{\theequation}{\thesection.\arabic{equation}}
\[
Department\; of\; Mathematics\; and\; Statistics
\]
\[
Wichita\; State\; University
\]
\[
Wichita,\; KS\; 67260-0033,\; U.S.A.
\]
\[
e-mail:\;\; victor.isakov@wichita.edu
\]

\abstract{ We study the problem of reconstruction of special special time dependent local volatility from market prices of options
with different strikes at two expiration times. For a general diffusion process we
apply the linearization technique and we conclude that
the option price can be obtained as the sum of the Black-Scholes
formula and of an operator $W$ which is linear
in perturbation of volatility. We further simplify the linearized inverse problem
and obtain unique solvability result in basic functional spaces. By using the Laplace
transform in time we simplify the kernels of integral operators
for $W$ and we obtain uniqueness and stability results for for volatility under natural condition of smallness of the spacial interval where one prescribes the (market) data. We propose a numerical algorithm based on our analysis of the linearized problem.}

\section{Introduction. Basic results}

The Black-Scholes formula \cite{BS} is an efficient 
method to price financial derivatives under the assumption
that the stock price is log-normally distributed. However, 
prices of options with different strikes generated by the Black-
Scholes formula differ from observed market prices \cite{H}. One way
to reconcile the differences is to replace the constant volatility by a more volatility
that might depnd on time and the underlying asset. This approach is very popular
in applications.

It is usually difficult to achieve a unique and stable fitting of the model
to actual market prices. While in the constant  volatility can be found very efficiently
 and uniquely from one option price, the space and time dependent volatility
function must be  restored from collection of simulteneous
option quotes with different strikes and expiration times. Even though many
numerical algorithms for this inverse problem have already 
been published \cite{AFHS}, \cite{BI1}, \cite{LO}, their convergence properties have not
been satisfactory for practitioners. In particular, since these algorithms use minimization of a non convex regularized misfit functional, there are well known difficulties with avoiding local minima, and their convergence is not guaranteed.
 We refer to \cite{BI2} for a survey of available theory and of numerical
methods. In the time dependent case a relation between implied and local volatilities was
discovered \cite{BBF}, \cite{BBF1} and a convergent numerical algorithm was designed.  In many practical situations the interval $\omega^*$ with strike prices is relatively small
and values of the volatility coefficient outside $\omega^*$ do not influence
option prices inside $\omega^*$ very much due to very fast diffusion. This
suggest a possibility of linearization. One of these linearizations was suggested in 
\cite{BJ} where linearized inverse problems are still complicated and no numerical
algorithms were discussed. In addition, expiration times are sparse, so it is hard to expect a detailed recovery of general time dependence of volatility.
In \cite{BIV} we assumed that volatility coefficient only depends on stock price, considered a linearized version of the inverse problem, obtained a simple convenient representation of this linearization, and desinged and tested a very fast and reliable numerical algorithm based on this linearization.  

In this paper we attempt to recover volatility which linearly depends on time while
the additional data are given for two expiration times. So we are making use of realistic
data to predict volatility in near future which seems to be most important for market
decisions. Below we outline the results in more detail.

  For any stock price, $0<s<\infty$, and time , $0<t<T$, the price $u$ for an option expiring at time $T$ satisfies the following partial differential equation
\begin{equation}
\frac{\partial u}{\partial t}+\frac{1}{2}s^2\sigma^2(s,t)
\frac{\partial^2 u}{\partial s^2}+
s\mu \frac {\partial u}{\partial s}-ru=0.
\label{BS}
\end{equation}
Here, $\sigma(s,t)$ is the volatility coefficient that satisfies
$0<m<\sigma(s,t)<M<\infty$ and is assumed to belong to
the H\"older space $C^{\lambda}(\bar\omega^*), 0<\lambda<1$ for some 
interval $\omega^*\subset (0,+\infty)$ and outside this interval, and $\mu$  
and $r$ are, respectively, the risk-neutral drift and the 
risk-free interest rate assumed to be constants.
 The backward in time parabolic equation (\ref{BS}) is 
augmented by the final condition specified by the pay off
of the call option with the strike price $K$ 
\begin{equation}
         u(s,T)=(s-K)^+=max(0,s-K), \; 0<s
\label{final}
\end{equation}
 It is known (e.g. see \cite{BI2} for H\"older $\sigma$) that 
there is a unique solution $u$ to (\ref{BS}),(\ref{final})
which belongs to $C^2((\bar\omega^*)\times (0,T]$ and
to $C((0,\infty)\times [0,T])$ and  satisfies the bound
$|u(s,t)|<C(s+1)$. The inverse problem of option pricing
seeks for $\sigma$ given
\begin{equation}
    u(s^*,t^*;K,T)=u^*(K,T),\; K\in\omega^*,\;T\in I^* .
\label{additionaldata}
\end{equation}
Here $s^*$ is market price of the stock at time $t^*$, 
 $u^*(K)$ denote market price of options with 
different strikes $K$ for a given expiry $T$. We
suppose that  $\omega^*$ contains $s^*$ and is small, and
attempt to recover volatility in the same interval.
 In this paper we assume that  $I^*=\{T_1, T_2\}$
 for some $0<T_1<T_2$ we are trying to recover the volatility
 coefficient
\begin{equation}
\frac{1}{2} \sigma^2(s)=\frac{1}{2}\sigma_0^2+f^*_0(s)+f^*_1(s)t
\label{perturbation}
\end{equation}
where $f^*_0, f^*_1$ are small $C(\bar\omega^*)$-perturbation of constant
$\sigma_0^2$ and $f^*_0=f^*_1=0$ outside $\omega_*$. 
Constant $\sigma_0$ is the implied volatility defined as the (unique) solution
to the Black-Scholes equation
\begin{equation}
\label{BSE}
s^*(N(\frac{1}{2}\sigma_0\sqrt{T_1-t^*})-
N(\frac{1}{2}(-\sigma_0\sqrt{T_1-t^*})))+u^8(s^*,T_1)
\end{equation}
where the normal distribution function
$$
N(x)=\frac{1}{\sqrt{2\pi}}\int_{-\infty}^x e^{-\frac{\theta^2}{2}} d\theta.
$$

Then the option price can be calculated (up to
a quadraticaly small error) as the sum of the Black-Scholes 
formula with volatility $\sigma_0$
and of a certain linear operator $U((f^*_0, f^*_1))$.

To obtain our results we will use that the option
premium $u(.,.;K,T)$ satisfies the equation dual to the
Black-Scholes equation (\ref{BS}) with respect
to the strike price $K$ and expiry time $T$:
\begin{equation}
\label{BS*}
\frac{\partial u}{\partial T}-
\frac{1}{2}K^2\sigma ^2(K, T) 
\frac{\partial ^2 u}{\partial K^2}+
\mu K\frac{\partial u}{\partial K}+( r-\mu) u=0  
\end{equation}
The equation (\ref{BS*}) was found by Dupire \cite{D} and rigorously justified, for example, in \cite{BI2}.

In the section 2 we use the standard linearization procedure
and derive a partial differential equation for the principal linear 
(with respect to $(f_0^*, f_1^*)$ part $V$ of $u$ in new logarithmic variables. 
After another  substitution this equation is reduced to the heat equation
with the right-hand side linear with respect to $(f^*_0, f^*_1)$.
We will drop $*$ in new variables.

Due to difficulties with the linearized inverse problem in section 3 we propose its  simplified version. This simplification is a good approximation in the most practically important case of small $\omega^*$. For the simplified inverse problem by using the Fourier transformation with respect to $s$  we obtain uniqueness, stability, and existence results in Sobolev spaces when the data are given for all strikes at two expire times.
Due to strong diffusion the data (\ref{additionaldata}) are small outside $\omega^*$, so 
the original data can be extended onto $(0,+\infty)\setminus\omega^*$ as some $C^2$-small function.
 
In section 4 we study the linearized inverse option pricing problem.
We apply the Laplace transform to explicitly evaluate the time integrals
in the standard integral representation of the solution of the direct option pricing problem. As a result we arrive in the following system of linear integral equations
for $(f_0, f_1)$ (obtained from $(f_0^*,f_1^*)$ via a simple transformation): 
 \begin{equation}
\label{integralequation}
\int_{\omega}(K_{0j}(x,y)f_0(y)+K_{1j}(x,y)f_1(y))dy=F_j(x),\;j=1,2,\; x\in\omega
\end{equation}
with the kernels 
$$
K_{0j}(x,y)=\frac{s^*}{\sigma^2_0 \sqrt{\pi}}
\int_{\frac{|x-y|+|y|}{ \sqrt{2\tau_j}\sigma_0}} e^{-\theta^2} d\theta,
$$ 
$$
K_{1j}(x,y)=\frac{s^*}{\sigma^2_0 \sqrt{\pi}}
(\frac{\sqrt{\tau_j}}{\sqrt{2}\sigma_0}(|y|-|x-y|)
e^{-\frac{(|x-y|+|y|)^2}{\sigma^2_0 2\tau_j}}+
$$
$$
(\frac{x^2-2xy}{\sigma_0^2}+\tau_j)
\int_{\frac{|x-y|+|y|}{ \sqrt{2\tau_j}\sigma_0}} e^{-\theta^2} d\theta,\tau_j=T_j-t^*, j=1,2,
$$
given by the error function and the right-hand side $F_j(x)$ computed from the market data
(\ref{additionaldata}).
Using this explicit formulae we prove that the system of integral equations
(\ref{integralequation}) has not more than one
solution under a certain constraint on
$\omega$ which is determined by the range of
strikes for which the data (\ref{additionaldata})
are given. It follows also a  that for small change in
$F_j$ that can be viewed as a market bid-ask spread or any
other noise in the market data (\ref{additionaldata}),
solutions changes remain small as well. In other words,
the linearized version of the volatility reconstruction is stable.

Numerical results obtained by other authors and by author 
of this paper using different methodologies (\cite{AFHS},\cite{BI1}, \cite{BIV},
\cite{EE}, \cite{LO}) exhibit similar reconstruction properties. In (\cite{AFHS}, \cite{LO}) they are looking for the volatility $\sigma=\sigma(t,s)$ from the data given for several $T$. We believe that values of the volatility for $t^*<t<T$ ( especially in the near future) are of most interest. But in the time dependent case the available
data are not sufficient to identify it. That is why we think one should
assume special time dependence and utilize other analytic features of this inverse problem. In \cite{BIV} we considered time independent $\sigma$ and use the data only at 
$T_1$.

Our current contribution is twofold.

First, we demonstrated uniqueness, stability, and existence of solution for the approximate linearized inverse problem. Furthemore,  we are able to explicitly describe (Theorem 4.1) the domain where the stable reconstruction is ensured for the linearized inverse problem. The closest result in this direction is obtained in 
\cite{BI2},  where the linearization was replaced by an approximate linearization, 
and in \cite{BIV} where only $f_=0$ and  $I^*=\{T_1\}$. 
 As in \cite{BIV}, the answer is obtained in terms of the (more complicated) ratios of the spatial variable  to the square root of time to maturity. It is very much consistent with trader's intuition.

Our approach suggests the actual reconstruction methodology. 
Algorithms proposed in \cite{AFHS}, \cite{BJ}, \cite{LO} attempt to use 
standard methods (regularized best $L^2$- fit) to numerical solution of 
inverse problems. They involve very heavy computations with poor convergence 
properties which is typical for inverse parabolic problems. The numerical
algorithm in \cite{BI1} is based on some features of the inverse option
pricing problem. In fact, it is a second order correction of the
Black-Scholes formula. Unfortunately, that this correction lacks a rigorous 
justification and numerically it is not very efficient. In (\cite{BJ})
they used perturbation methods which look reasonable when volatility
is close to a constant, but the corresponding 
linear integral equation remains complicated as well as  its proposed 
numerical inversion.  In \cite{BIV} we were able to transform the linearized inverse problem into a one-dimensional integral equation with an explicit kernel.
This led to a very fast, simple, and stable numerical algorithm for the 
reconstruction of the local volatility. 
In this paper we transform the linearized inverse option pricing problem
into the system of linear integral equations (\ref{integralequation}). 
We already started numerical solution of this system.

\section{Linearization at constant volatility}

The substitution
$$
y=\ln ( \frac{K}{s^*}) ,\quad \tau =T-t^*,
$$
\begin{equation}
a(y,\tau)=\sigma(s^*e^y,\tau+t^*),\;U(y,\tau)=u(s^*e^y,\tau+t^*)  
\label{substitution}
\end{equation}
transforms the equation (\ref{BS*}) and the initial data (\ref{final}) into
$$
\frac{\partial U}{\partial \tau}-
\frac{1}{2}a^2 ( y,\tau) \frac{\partial^2 U}{\partial y^2}+
( \frac{1}{2} a^2( y,\tau) +
\mu ) \frac{\partial U}{\partial y}+( r -\mu ) U =0,
$$
\begin{equation} 
 \label{BS*y} 
U( y,0)  = s^{*}( 1-e^y)^{+},\; y\in \mathbb{R}. 
\end{equation}
while the additional (market) data become
\begin{equation}
U( y,\tau_j ) =U_j(y),\quad y\in \omega, \;j=1,2.  
\label{final*}
\end{equation}
Here $\omega$ is the transformed interval $\omega^*$ 
($\omega^*$ in $y$- variables (\ref{substitution})) and
$U_j(y)=u^*(s^* e^y, T_j)$. Observe that $\tau_j=T-t_j^*$.

The equations (\ref{BS*y}) and (\ref{final*}) for functions $U(\tau, y), a(y)$
form the so-called inverse parabolic problem with the final overdetermination.
The known uniqueness conditions for this problem \cite{Is1}, section 6.2, 
\cite{Is2}, section 9.2, are not satisfied in our particular situation, 
and we are not aware of any uniqueness result.

To derive the linearized inverse problem we observe that due to the
assumption (\ref{perturbation}),
$$
\frac{1}{2} a^2( y,\tau) =\frac{1}{2}\sigma _0^2+f_0(y)+\tau f_1(y),   
$$
where 
$$
f_0(y)=f_0^*(s^*e^y)+t^*f_1^*(s^*e^y), f_1(y)=f_1^*(s^*e^y)
$$ 
are $C(\bar\omega)$-small and equal to zero outside $\omega$. So
\begin{equation}
U=V_0+V+v.  \label{splitting}
\end{equation}
Here $V_0$ solves (\ref{BS*y}) with $a=\sigma _0$ and 
$v$ is quadratically small with respect to $(f_0,f_1)$, while the
principal linear term $V$ satisfies the equations
$$
\frac{\partial V}{\partial \tau }-
\frac 12\sigma _0^2\frac{\partial ^2V}{\partial y^2}+( \frac{\sigma _0^2}2+\mu ) \frac{\partial V}{\partial y}+( r -\mu ) V = \alpha_0 (f_0(y)+\tau f_1(y)),
$$
\begin{equation}  
V( y,0 )  =0 ,\quad y\in \mathbb{R}, 
\label{BSlinear} 
\end{equation}
where 
$$
\alpha _0( y,\tau )  =\frac {s^*}{\sqrt{2 \pi \tau} \sigma_0}
e^{-\frac{y^2}{2\tau \sigma _0^2}+cy+d\tau }, 
$$
$$
c =\frac 12+\frac \mu {\sigma _0^2},\quad 
d=-\frac 1{2\sigma _0^2}( \frac{\sigma _0^2}2+\mu ) ^2+\mu -r 
$$
and the additional final data
\begin{equation}
\label{additionaldataV}
V(y,\tau_j)=V_j(y), j=1,2, y\in\omega,
\end{equation}
where $V_j(y)=U_j(y)-V_0(y,\tau_j)$.
One can completely justify this linearization by using standard theory of direct
parabolic boundary  value problems \cite{F}, \cite{LSU}, as it was done in \cite{BI2},
 pp. R 103-104, for the inverse option pricing problem or in
\cite{Is2}, section 4.5, for some elliptic inverse problems. 

The new substitution
\begin{equation}
V=e^{cy+d\tau }W 
\label{newsubstitution}
\end{equation}
simplifies (\ref{BSlinear}) to
$$
\frac{\partial W}{\partial \tau }-
\frac 12\sigma_0^2\frac{\partial ^2W}{\partial y^2} =
\alpha (f_0(y)+\tau f_1(y)),\quad 0< \tau <\tau ^*, \quad y \in {\mathbb R},
$$ 
\begin{equation}
  \label{heat} 
W( y,0)  =0,\quad y\in {\mathbb R},
\end{equation}  
where
$$
\alpha ( \tau ,y)  =
\frac{s^*}{\sqrt{2\pi \tau }\sigma _0}
e^{-\frac{y^2}{2\tau \sigma _0^2}},
$$
with the additional final data
\begin{equation}
\label{additionaldataW}
W(y,\tau_j)=W_j(y),\quad y\in\omega, j=1,2. 
\end{equation}
where
$$
W_j(y)=e^{-cy-d\tau_j}V_j(y).
$$

Due to analytic difficulties, using some features of our inverse problem we will further simplify (\ref{BSlinear}) to the approximate linearized inverse problem by replacing $\alpha$ with the best approximating $y$ independent function $\alpha_1$:
$$
\frac{\partial w}{\partial \tau }-
\frac 12\sigma _0^2\frac{\partial ^2W}{\partial y^2} 
=\alpha_1 (f_0(y)+\tau f_1(y)),\quad 0< \tau <\tau^*, \quad y \in {\mathbb R}, 
$$
\begin{equation}
\label{heatA}
w( y,0)  =0,\quad y\in {\mathbb R}
\end{equation}
where
$$
\alpha_1( \tau ,y)  =
\frac{s^*}{\sqrt{2\pi \tau }\sigma _0}
$$
with the additional final data
\begin{equation}
\label{additionaldataA}
W(y,\tau_j)=W_j(y),\quad y\in\mathbb R, j=1,2. 
\end{equation}

In particular for computational purposes we replace $\mathbb R$ by a finite interval $\Omega$ containing $\bar\omega$:
$$
\frac{\partial w}{\partial \tau }-
\frac 12\sigma _0^2\frac{\partial ^2W}{\partial y^2} 
=\alpha_1 (f_0(y)+\tau f_1(y)),\quad 0< \tau <\tau^*, \quad y \in \Omega, 
$$
\begin{equation}
\label{heatA1}
w( y,0)  =0,\quad y\in \Omega
\end{equation}
with the additional final data
\begin{equation}
\label{additionaldataA1}
W(y,\tau_j)=W_j(y),\quad y\in\Omega, j=1,2. 
\end{equation}

As in \cite{BI2}, Theorem 4, one can show that uniqueness and stability for the approximate inverse problem imply uniqueness and stability for the linearized inverse problem when the interval $\omega$ is small.

\section{The approximate linearized inverse problem}

In this section we prove the following unique solvability result

\begin{theorem}

Let $W_j\in H^{(2)}(\mathbb{R}), j=1,2$. 

Then there is a unique solution 
$(f_0,f_1)\in H^{(0)}(\mathbb{R})\times H^{(0)}(\mathbb{R})$ to the approximate linearized inverse problem (\ref{heatA}), (\ref{additionaldataA}).

Moreover there is a constant $C$ depending only on $\tau_1,\tau_2, \sigma_0$ such that a solution $(f_0,f_1)$ to the approximate linearized inverse problem satisfies the bound
\begin{equation}
\label{boundA}
\|f_0\|_{(0)}(\mathbb{R})+\|f_1\|_{(0)}(\mathbb{R})\leq
C(\|W_1\|_{(2)}(\mathbb{R})+\|W_2\|_{(2)}(\mathbb{R}))
\end{equation}
\end{theorem}

We use  the Fourier transform 
$$
{\cal F}f(\xi, \tau)=(2\pi)^{-\frac{1}{2}}\int f(y,\tau)e^{-i\xi y} dy
$$ 
of a function $f(y,\tau)$ with respect to $y$ and the standard norm
$$
\|f\|_{(l)}(\mathbb{R})=(\int(1+\xi^2)^l|{\cal F} f(\xi)|^2 d\xi)^{\frac{1}{2}}.
$$
of a function $f$ in the Sobolev space $H^{(l})(\mathbb{R})$.

{\bf Proof:}

Applying the Fourier transform to the initial value problem (\ref{heatA}) we yield
$$
\partial_{\tau} ({\cal F} W)+ \Xi^2  {\cal F} W = 
S(\tau^{-\frac{1}{2}} {\cal F} f_0+ \tau^{\frac{1}{2}}{\cal F} f_1),\; 0<\tau<\tau^*_j;
\quad {\cal F} W( ,0)=0,
$$
where we let
$$
S=\frac{s^*}{\sigma_0 \sqrt{2\pi}}, \; \Xi = \frac{1}{\sqrt{2}} \sigma_0 \xi.
$$
Solving the initial value problem for the ordinary differential equation we obtain
$$
{\cal F} W (\xi,\tau)=
$$
\begin{equation}
\label{W}
S(\int_0^{\tau}\theta^{-\frac{1}{2}}e^{\Xi^2(\theta-\tau)} d\theta) {\cal F} f_0(\xi)+
S(\int_0^{\tau}\theta^{\frac{1}{2}}e^{\Xi^2(\theta-\tau)}d\theta){\cal F} f_1(\xi).
\end{equation}
Hence  the inverse problem (\ref{heatA}),(\ref{additionaldataA})
is equivalent to the system of two integral equations
$$
(\int_0^{\tau_1}\theta^{-\frac{1}{2}}e^{\Xi^2(\theta-\tau_1)}d\theta){\cal F} f_0(\xi)+
(\int_0^{\tau_1}\theta^{\frac{1}{2}}e^{\Xi^2(\theta-\tau_1)}d\theta){\cal F} f_1(\xi) =
$$
$$
 S^{-1}{\cal F} W_1(\xi),
$$
$$
(\int_0^{\tau_2}\theta^{-\frac{1}{2}}e^{\Xi^2(\theta-\tau_2)}d\theta){\cal F} f_0(\xi)+
(\int_0^{\tau_2}\theta^{\frac{1}{2}}e^{\Xi^2(\theta-\tau_2)}d\theta){\cal F} f_1(\xi) =
$$
\begin{equation}
\label{systemA} 
S^{-1} {\cal F} W_2(\xi).
\end{equation}
Solving this linear algebraic system we obtain
$$
{\cal F} f_0(\xi)=
$$
$$
d^{-1}(\xi)((\int_0^{\tau_2}\theta^{\frac{1}{2}}
e^{\Xi^2(\theta-\tau_2)}d\theta){\cal F} W_1(\xi)-
(\int_0^{\tau_1}\theta^{\frac{1}{2}}e^{\Xi^2(\theta-\tau_1)}d\theta {\cal F} W_2(\xi),
$$
$$
{\cal F} f_1(\xi)=
$$
\begin{equation}
\label{f}
d^{-1}(\xi)((\int_0^{\tau_1}
\theta^{-\frac{1}{2}} e^{\Xi^2(\theta-\tau_1)}d\theta) {\cal F} W_2(\xi)-
(\int_0^{\tau_2}\theta^{-\frac{1}{2}}e^{\Xi^2(\theta-\tau_2)}d\theta {\cal F} W_1(\xi)
\end{equation}
where
$$
d(\xi)=S 
(\int_0^{\tau_1}\theta^{-\frac{1}{2}}e^{\Xi^2(\theta-\tau_1)}d\theta 
\int_0^{\tau_2}\theta^{\frac{1}{2}}e^{\Xi^2(\theta-\tau_1)}d\theta-
$$
$$
\int_0^{\tau_1}s^{\frac{1}{2}}e^{\Xi^2(\theta-\tau_1)}d\theta \int_0^{\tau_2}\theta^{-\frac{1}{2}}e^{\Xi^2 (\theta-\tau_2)}d\theta).
$$

Integrating by parts,
$$
\int_0^{\tau}\theta^{\frac{1}{2}}e^{\Xi^2(\theta-\tau)}d\theta = 
\frac{\tau^{\frac{1}{2}}}{\Xi^2}-
\frac{1}{2\Xi^2}\int_0^{\tau}\theta^{-\frac{1}{2}}e^{\Xi^2(\theta-\tau)}d\theta.
$$
By elementary substitution $\theta=\tau_j\rho$,
$$
\int_0^{\tau_j}\theta^{-\frac{1}{2}}e^{\Xi^2(\theta-\tau_j)}d\theta
=
\tau_j^{\frac{1}{2}}\int_0^{1}\rho^{-\frac{1}{2}}e^{\Xi^2\tau_j(\rho-1)}d\rho
$$
and therefore
\begin{equation}
\label{d}
d(\xi)=S \frac{\sqrt{\tau_1\tau_2}}{\Xi^2}
(\int_0^1\rho^{-\frac{1}{2}}e^{\Xi^2\tau_1(\rho-1)}d\rho -
\int_0^1\rho^{-\frac{1}{2}}e^{\Xi^2\tau_2(\rho-1)}d\rho).
\end{equation}

Using the Taylor expansion for the exponential function and the definition of $\Xi=
\frac{1}{\sqrt{2}}\sigma_0\xi$ we have
$$
e^{\Xi^2\tau_j(\rho-1)}=1+\Xi^2\tau_j(\rho-1)+O(\xi^4),
$$
where $|O(\xi)|\leq C |\xi|$, when $|\xi|<1$. So calculating the integral of
$\rho^{\frac{1}{2}} - \rho^{-\frac{1}{2}}$ over $(0,1)$ we yield
\begin{equation}
\label{d1}
d(\xi)=S \frac{4}{3} \sqrt{\tau_1 \tau_2}(\tau_2-\tau_1)+O(\xi^2) \;
\mbox{provided} \;|\xi|<1.
\end{equation}

Splitting the integration interval and integrating by parts again in the first integral,
we will have
$$
  \int_0^1\rho^{-\frac{1}{2}}e^{\Xi^2 \tau_j(\rho-1)}d\rho =
  \int_{0.5}^1\rho^{-\frac{1}{2}}e^{\Xi^2 \tau_j(\rho-1)}d\rho +
  \int_0^{0.5}\rho^{-\frac{1}{2}}e^{\Xi^2 \tau_j(\rho-1)}d\rho  =
$$
$$
\frac{1}{\tau_j\Xi^2}+
\frac{1}{2\tau_j\Xi^2}\int_{0.5}^1\rho^{-\frac{3}{2}}
e^{\Xi^2\tau_j(\rho-1)}d\rho + O(\xi^{-4}) =
\frac{1}{\tau_j\Xi^2}+ O(\xi^{-4}).
$$
when $1\leq |\xi|$, if we repat the integration by parts. Hence
\begin{equation}
\label{d2}
d(\xi)=S \frac{\tau_2-\tau_1}{\sqrt{\tau_1 \tau_2}\Xi^4}+O(\xi^{-6}) \;
\mbox{provided} \;|\xi|\geq 1.
\end{equation}

Using that $\Xi= \frac{1}{\sqrt{2}}\sigma_0\xi$, from 
 \eqref{d1}  \eqref{d2} we derive that
\begin{equation}
\label{d3}
C^{-1}(1+\xi^2)^{-2}\leq d(\xi)\leq C(1+\xi^2)^{-2}.
\end{equation}

Similarly to the derivation of \eqref{d2} we obtain
$$
\int_0^{\tau_j}\theta^{-\frac{1}{2}}e^{\Xi^2(\theta-\tau_j)}d\theta=
\frac{1}{\sqrt{\tau_j}\Xi^2}+O(\xi^{-4}),
$$
$$
\int_0^{\tau_j}\theta^{\frac{1}{2}}e^{\Xi^2(\theta-\tau_j)}d\theta=
\sqrt{\tau_j}\Xi^{-2}+O(\xi^{-4}),
$$
if $|\xi|>1$, and therefore,
$$
|\int_0^{\tau_j}\theta^{-\frac{1}{2}}e^{\Xi^2(\theta-\tau_j)}d\theta|\leq C(1+\xi^2)^{-1},
$$
\begin{equation}
\label{integral}
|\int_0^{\tau_j}\theta^{\frac{1}{2}}e^{\Xi^2(\theta-\tau_j)}d\theta|\leq C(1+\xi^2)^{-1}.
\end{equation}

From \eqref{f}, \eqref{d}, \eqref{d3}, and \eqref{integral} we conclude that
$$
|{\cal F} f_j(\xi)|^2\leq C(1+\xi^2)^2(|{\cal F} W_1|^2(\xi)+|{\cal F} W_2|^2(\xi)),
\; j=0,1.
$$
Using the definition of Sobolev norms via the Fourier transform we obtain the bound
\eqref{boundA}. This bound implies uniqueness and stability. 

The formulae \eqref{f} give a solution explicitly for smooth compactly supported
data $(W_1,W_2)$, combined with the density of such data in 
$H^{(2)}({\mathbb R})\times H^{(2)}({\mathbb R})$ and the stability estimate 
\eqref{boundA} it implies  existence of solution in the needed Sobolev space.

The proof is complete.

\begin{theorem}

Let $0<B, \Omega=(-B, B), \;W_j\in H^{(2)}(\Omega),  W_j(-B)= W_j (B)=0,\; j=1,2$. 

Then there is a unique solution 
$(f_0,f_1)\in H^{(0)}(\Omega)\times H^{(0)}(\Omega)$ to the approximate linearized inverse problem (\ref{heatA1}), (\ref{additionaldataA1}).

Moreover there is a constant $C$ depending only on $\tau_1,\tau_2, B , \sigma_0$ such that a solution $(f_0,f_1)$ to the approximate linearized inverse problem satisfies the bound
$$
\|f_0\|_{(0)}(\Omega)+\|f_1\|_{(0)}(\Omega)\leq
C(\|W_1\|_{(2)}(\Omega)+\|W_2\|_{(2)}(\Omega))
$$
\end{theorem}

The proof is similar to the above proof of Theorem 3.1 if instead of the Fourier transform we use orthonormal (eigenfunction) series
$$
f(y,\tau)=\sum_{n=1}^{\infty} f_n(\tau) \frac{1}{\sqrt{B}} sin \frac{\pi n(y+B)}{2B},\;
f_n(\tau)= \int_{-B}^B  f(y,\tau)\frac{1}{\sqrt{B}} sin \frac{\pi n(y+B)}{2B}.
$$
It suffices to replace $\xi$ by $n$.

By elementary means ( using sharp bounds in Taylor formula ) one can explicitly evaluate constants $C$. We did not do it to simplify the exposition.

\section{Uniqueness and stability for the linearized inverse problem}

The purpose of this section is to prove the following

\begin{theorem}
Let $\omega=\left( -b,b\right)$ and $f_0=f_1=0$ outside $\omega$. 

If 
$$
\frac{\tau_1^2+\tau_2^2+\sqrt{\tau_1\tau_2}(\tau_1+\tau_2)}{\tau_1\tau_2(\tau_2-\tau_1)}
\frac{3b^2}{2\sigma_0^2}<1.
$$
\begin{equation}
\frac{(\sqrt{\frac{\tau_1}{\tau_2}}+\sqrt{\frac{\tau_2}{\tau_1}}+2)
\frac{b^2}{\sigma_0^2}+
2\sqrt{2\pi}(\sqrt{\tau_1}+\sqrt{\tau_2})\frac{b}{\sigma_0}}
{2(\tau_2-\tau_1)}< 1,
\label{uniqueness}
\end{equation}
 then a solution 
$(f_0,f_1)\in L^\infty(\omega)\times L^\infty(\omega)$ to the linearized inverse
option pricing problem \eqref{BSlinear}, \eqref{additionaldataV} is unique.

Moreover, there is a constant $C$ depending only on $\sigma_0, \tau_1, \tau_2, b$,
such that
\begin{equation}
\label{boundB}
\|f_0\|_{\infty}(\mathbb{R})+\|f_1\|_{\infty}(\mathbb{R})\leq
C(\|W_1"\|_{\infty}(\mathbb{R})+\|W_2"\|_{\infty}(\mathbb{R})).
\end{equation}
\end{theorem}

We remind that $\|f\|_{\infty}(\omega)$ is
essential supremum of $|f|$ over $\omega$. In particular,
when $f$ is continuous on $\overline{\omega}$ is it just
$max|f(x)|$ over $x\in \overline{\omega}$.

 From  numerical experiments
in \cite{LO} and in \cite{BIV}  we can see that
values of $\sigma$ outside $\omega$ are not essential, due to
a very fast decay of the Gaussian kernel $\alpha$ in $y$.

\begin{lemma}

For a solution $W(\cdot;f_0,f_1)$ of the Cauchy problem \eqref{heat} we have
\begin{equation}
W(x,\tau; f_0,f_1) =\int_{\omega}(K_0( x,y;\tau ) f_0( y)+K_1( x,y;\tau ) f_1( y)) dy,\quad
x\in \omega,  
\label{Wintegral}
\end{equation}
where
$$
K_0( x,y;\tau) =S_*\int_{\frac{|x-y| +| y| }{\sigma _0\sqrt{2\tau}}}^{\infty} 
e^{-\tau ^2}d\tau, \; S_*=\frac{s^*}{\sigma^2_0\sqrt{\pi}}, 
$$
$$
K_1( x,y;\tau) =S_*(\frac{\sqrt{\tau}}{\sqrt{2}\sigma_0}(|y|-|x-y|)
e^{-\frac{(|x-y| +| y|)^2 }{2\tau \sigma_0^2}}+
$$
$$
(\frac{x^2-2xy}{\sigma_0^2}+\tau)
\int_{\frac{|x-y| +| y| }{\sigma _0\sqrt{2\tau}}}^{\infty} 
e^{-\tau ^2}d\tau). 
$$
\end{lemma}

{\bf Proof: } 

When $f_1=0$ this result was obtained in \cite{BIV}. So by using the linearity with respect to $(f_0, f_1)$ we can assume that $f_0=0$.

The well-known representation \cite{F} of the solution to the 
Cauchy problem (\ref{heat}) for the heat equation yields
$$
W( x,\tau)  = \int_{\mathbb R} K( x,y;\tau) f_1(y) dy,
$$
where
$$
K( x,y;\tau)  = \frac{s^*}{2\pi\sigma_0^2}
\int_0^{\tau} \frac {1}{\sqrt{( \tau-\theta ) }}
e^{-\frac{| x-y| ^2}{2\sigma_0^2( \tau -\theta ) }}
\frac {1}{\sqrt{\theta }}
\theta e^{-\frac{| y| ^2}{2\sigma _0^2\theta }}d\theta. 
$$
 
As in \cite{BIV} we will simplify $K( x,y;\tau) $ by using the Laplace
transform $\Phi ( p) ={\cal L} (\phi)( p ) $ of  $\phi ( \tau ) $ with respect to 
$\tau $. Since the Laplace transform of the convolution is the product of Laplace transforms of convoluted functions, we have
$$
{\cal L} K ( x,y; )( p)  = \frac {s^*}{2\pi \sigma _0^2}
{\cal L}( \frac 1{\sqrt{\tau }}e^{-\frac{| x-y| ^2}{2\sigma _0^2\tau }}) 
{\cal L}( \sqrt{\tau }e^{-\frac{| y| ^2}{2\sigma _0^2\tau }}) =
$$
$$
\frac {s^*}{2\pi \sigma _0^2}\sqrt{\frac {\pi}{ p}}
e^{-\frac{\sqrt{2}|x-y| }{\sigma _0}\sqrt{p}}
\frac{\sqrt{ \pi}}{ 2}
e^{-\frac{\sqrt{2}|y| }{\sigma _0}\sqrt{p}}(p^{-\frac{3}{2}}+
\sqrt{2}\frac{|y|}{\sigma_0 p})=
$$
\begin{equation}
\label{LK}
=\frac {s^*}{2\sigma _0^2}
(p^{-2} e^{-\frac{\sqrt{2}(| x-y|+| y|)}{\sigma _0}\sqrt{p}}+
\frac{\sqrt{2}|y|}{\sigma_0}p^{-\frac{3}{2}} 
e^{-\frac{\sqrt{2}(| x-y|+| y|)}{\sigma _0}\sqrt{p}}),
\end{equation}
where we used the formula for the Laplace transform of
$\tau^{-\frac{1}{2}} e^{-\frac{\beta}{\tau}}$ \cite{LS},
p.526, formula 16, and that ${\cal L}(\tau\phi(\tau))=-\frac{d}{dp}{\cal L}(\phi(\tau))$,
\cite{LS}, p.499.
  
Applying the formula for the inverse Laplace transform of the 
function  $\frac{1}{p}e^{-\gamma \sqrt{p}}$ \cite{LS}, 
p. 528, formula 42, and using the relation between Laplace transforms of a function
and its integral 
$$
{\cal L}^{-1}(p^{-1}\Phi(p))=\int_0^{\tau}\phi(\theta) d\theta,
$$
\cite{LS}, p. 499, and \cite{LS}, p. 526, formula 16, we conclude that
$$
{\cal L}^{-1}(\frac{e^{-\gamma\sqrt{p}}}{p^2})=
\frac{2}{\sqrt{\pi}}
\int_0^{\tau}(\int_{\frac{\gamma}{2\sqrt{\theta}}} e^{-\rho^2} d\rho) d\theta=
$$
$$
\frac{2}{\sqrt{\pi}}
\int_{\frac{\gamma}{2\sqrt{\tau}}}^{\infty}e^{-\rho^2}
(\int_{(\frac{\gamma}{2{\rho}})^2}^{\tau}1 d\theta) d\rho=
\frac{2}{\sqrt{\pi}}
\int_{\frac{\gamma}{2\sqrt{\tau}}}^{\infty}(\tau-(\frac{\gamma}{2\rho})^2 )e^{-\rho^2}
 d\rho=
$$
$$
\frac{2}{\sqrt{\pi}}
(\tau\int_{\frac{\gamma}{2\sqrt{\tau}}}^{\infty}e^{-\rho^2} d\rho -
\frac{\gamma^2 2\sqrt{\tau}}{4\gamma} e^{-\frac{\gamma^2}{4\tau}}+
\frac{\gamma^2}{2} \int_{\frac{\gamma}{2\sqrt{\tau}}}^{\infty}e^{-\rho^2} d\rho)=
$$
\begin{equation}
\label{L1}
\frac{2}{\sqrt{\pi}}
((\tau+\frac{\gamma^2}{2})\int_{\frac{\gamma}{2\sqrt{\tau}}}^{\infty}e^{-\rho^2} d\rho-
\frac{\gamma \sqrt{\tau}}{2} e^{-\frac{\gamma^2}{4\tau}})
\end{equation}
where we switched order of integration with respect to $\rho$ and $\theta$ and integrated by parts.

To find the inverse Laplace transform for \eqref{LK} we will again use
 the known above relation between Laplace transform of a functions and of its integral, to find
$$
{\cal L}^{-1}(p^{-\frac{3}{2}}e^{-\gamma\sqrt{p}})=
\int_0^{\tau}\frac{1}{\sqrt{\pi \theta}} e^{-\frac{\gamma^2}{4\theta}} d\theta =
\frac{\gamma}{\sqrt{\pi}}\int_{\frac{\gamma}{2\sqrt{\tau}}}^{\infty}\rho^{-2}) e^{-\rho^2} d\rho =
$$
\begin{equation}
\label{L2}
\frac{\gamma}{\sqrt{\pi}}\int_{\frac{\gamma}{2\sqrt{\tau}}}^{\infty}
\partial_{\rho}(-\rho^{-1}) e^{-\rho^2} d\rho =
\frac{2\gamma}{\sqrt{\pi}}(\frac{\sqrt{\tau}}{\gamma} e^{-\frac{\gamma^2}{4\tau}}-
\int_{\frac{\gamma}{2\sqrt{\tau}}}^{\infty} e^{-\rho^2} d\rho 
\end{equation}
where we utilized the substitution $\theta=\frac{\gamma^2}{4\rho^2}$ and integrated by parts.

By using \eqref{L1}, \eqref{L2} with $\gamma=\sqrt{2}\frac{|x-y|+|y|}{\sigma_0}$ we yield
$$
{\cal L}^{-1}(\frac {s^*}{2\sigma _0^2}
(p^{-2} e^{-\frac{\sqrt{2}(| x-y|+| y|)}{\sigma _0}\sqrt{p}}+
\frac{\sqrt{2}|y|}{\sigma_0}p^{-\frac{3}{2}} 
e^{-\frac{\sqrt{2}(| x-y|+| y|)}{\sigma _0}\sqrt{p}}))=
$$
 $$
\frac {s^*}{2\sigma _0^2}
(\frac{2}{\sqrt{\pi}}((\tau+\frac{(|x-y|+|y|)^2}{\sigma_0^2})
\int_{\frac{|x-y| +| y| }{\sigma _0\sqrt{2\tau}}}^{\infty} 
e^{-\rho ^2}d\rho-
$$
$$
\sqrt{2}\frac{| x-y|+| y|}{2\sigma _0}\sqrt{\tau} 
e^{-\frac{(| x-y|+| y|)^2}{2\tau\sigma_0^2}})+
$$
$$ 
\frac{\sqrt{2}| y|}{\sigma _0}(2\sqrt{\frac{\tau}{\pi}} 
e^{-\frac{(| x-y|+| y|)^2}{2\tau\sigma_0^2}}-
2\sqrt{2}\frac{|x-y|+|y|}{\sqrt{\pi}\sigma_0}
\int_{\frac{|x-y| +| y| }{\sigma _0\sqrt{2\tau}}}^{\infty} 
e^{-\rho ^2}d\rho))=
$$
$$
\frac {s^*}{2\sigma _0^2}(\frac{\sqrt{2\tau}}{\sigma_0\sqrt{\pi}}
(|y|-| x-y|) e^{-\frac{(| x-y|+| y|)^2}{2\tau\sigma_0^2}}+
\frac{2}{\sqrt{\pi}}(\frac{x^2-2xy}{\sigma_0^2}+\tau)
\int_{\frac{|x-y| +| y| }{\sigma _0\sqrt{2\tau}}}^{\infty} 
e^{-\rho ^2}d\rho).
$$
Now the formula of Lemma 4.2 for $K_1$ follows from \eqref{LK}.

\begin{lemma}
   The linearized inverse problem \eqref{heat}, \eqref{additionaldataW} implies the following Fredholm system of the linear integral equations
$$
f_0+\tau_1 f_1+ A_{11}f_0+A_{12} f_1= w_1,
$$
\begin{equation}
\label{system}
f_1+\tau_2 f_1+ A_{21}f_0+A_{22} f_1= w_2,
\end{equation}
where
$$
A_{j1}f_0(x)=-\frac{1}{2\tau_j\sigma_0^2}
\int_{\omega} (|x-y|+|y|)e^{-\frac{( |x-y| +| y|) ^2-x^2}{2\tau_j\sigma _0^2}}
f_0(y) dy, \;j=1,2,
$$
\begin{equation}
\label{A}
A_{j2}f_1(x)=- \int_{\omega}( \frac{|y|}{\sigma_0^2}
e^{-\frac{( |x-y| +| y|) ^2-x^2}{2\tau_j\sigma _0^2}}+
\sqrt{\frac{\tau_j}{2}}\frac{1}{\sigma_0}e^{\frac{x^2}{2\tau_j\sigma _0^2}} 
 \int_{\frac{|x-y|+|y|}{\sigma_0\sqrt{2\tau_1}}}^{\infty}e^{-\theta^2} d\theta) f_1(y) dy,
\end{equation}
and
$$
w_j(x)=-\frac{\sqrt{\pi\tau_j}\sigma_0^3}{\sqrt{2}s^*}
e^{\frac{x^2}{2\tau_j\sigma _0^2}} \partial_x^2W_j(x).
$$
\end{lemma}

{\bf Proof:}

The formula for $A_{j1}f_0$ was derived in \cite{BIV}.

Now we will consider $A_{j2}f_1$.

Using that $\partial_x|x-y|=sign(x-y)$ we obtain
$$
\partial_x \int_{\omega} K_1( x,y;\tau ) f_1( y) dy =
$$
$$
\frac{s^*}{2\sigma_0^2\sqrt{\pi}}\int_{\omega}(\frac{\sqrt{\tau}}{\sqrt{2}\sigma_0}
(-sign(x-y))e^{-\frac{(|x-y| +| y|)^2 }{2\tau \sigma_0^2}}+
$$
$$
\frac{\sqrt{\tau}}{\sqrt{2}\sigma_0}
(|y|-|x-y|)e^{-\frac{(|x-y| +| y|)^2 }{2\tau \sigma_0^2}}
(-\frac{|x-y|+|y|}{\tau\sigma_0^2}sign(x-y))+
$$
$$
(\frac{x^2-2xy}{\sigma_0^2}+\tau)
(-e^{-\frac{(|x-y| +| y|)^2 }{2\tau \sigma_0^2}}
\frac{sign(x-y)}{\sqrt{2\tau}\sigma_0})+
$$
$$
\frac{2x-2y}{\sigma_0^2}
\int_{\frac{|x-y| +| y| }{\sigma _0\sqrt{2\tau}}}^{\infty} 
e^{-\theta ^2}d\theta) f_1(y) dy = 
$$
$$
\frac{s^*}{\sigma_0^3\sqrt{\pi}}\int_{\omega}
(-\frac{\sqrt{\tau}}{\sqrt{2}}
sign(x-y)e^{-\frac{(|x-y| +| y|)^2 }{2\tau \sigma_0^2}}+
\frac{x-y}{\sigma_0}
\int_{\frac{|x-y| +| y| }{\sigma _0\sqrt{2\tau}}}^{\infty} 
e^{-\theta ^2}d\theta) f_1(y) dy . 
$$

Differentiating once more and using that $\partial_x sign(x-y)=2 \delta(x-y)$
we yield
$$
\partial^2_x \int_{\omega} K_1( x,y;\tau ) f_1( y) dy =
\frac{s^*}{\sigma_0^3\sqrt{\pi}}
(-\frac{\sqrt{\tau}}{\sqrt{2}}2 f_1(x)e^{-\frac{x^2 }{2\tau \sigma_0^2}}+
$$
$$
\int_{\omega}
(-\frac{\sqrt{\tau}}{\sqrt{2}} sign(x-y)
e^{-\frac{(|x-y| +| y|)^2 }{2\tau \sigma_0^2}}
(-\frac{|x-y|+|y|}{\tau\sigma_0^2} sign(x-y))+
$$
$$
\frac{1}{\sigma_0}
\int_{\frac{|x-y| +| y| }{\sigma _0\sqrt{2\tau}}}^{\infty} e^{-\theta ^2}d\theta
+  
\frac{x-y}{\sigma_0}
 e^{-\frac{(|x-y| +| y|)^2 }{2\tau\sigma_0^2}}
 (-\frac{1}{\sqrt{2\tau}\sigma_0} sign(x-y))) f_1(y) dy =  
$$
$$
\frac{s^*}{\sigma_0^3\sqrt{\pi}}
(-\sqrt{2\tau} f_1(x)e^{-\frac{x^2 }{2\tau \sigma_0^2}}+
$$
$$
\int_{\omega}
(\frac{|y|}{\tau\sigma_0^2} 
e^{-\frac{(|x-y| +| y|)^2 }{2\tau \sigma_0^2}}+
\frac{1}{\sigma_0}
\int_{\frac{|x-y| +| y| }{\sigma _0\sqrt{2\tau}}}^{\infty} e^{-\theta ^2}d\theta) 
f_1(y) dy) .  
$$

The proof is complete.

In our opinion, several cancellations in the the proof of Lemma 4.3 suggest that our approach is quite suitable for this problem

{\bf Proof of Theorem 4.1: }

Due to Lemma 4.3 to prove Theorem 4.1 it suffices to show uniqueness and stability
of solution $(f_0,f_1)$ of (\ref{system}). By linear algebra 
\eqref{system} is equivalent to the system
$$
f_0+\frac{1}{\tau_2-\tau_1}
((\tau_2A_{11}-\tau_1A_{21})f_0+(\tau_2A_{12}-\tau_1A_{22}) f_1)= 
\frac{1}{\tau_2-\tau_1}(\tau_2 w_1-\tau_1 w_2),
$$
\begin{equation}
\label{systemd}
f_1+\frac{1}{\tau_2-\tau_1}
((-A_{11}+A_{21})f_0+(-A_{12}+A_{22}) f_1)= 
\frac{1}{\tau_2-\tau_1}(- w_1 + w_2).
\end{equation}

Using \eqref{A},
$$
\|A_{j1}f_0\|\leq 
  \frac{1}{2\tau_j\sigma_0^2} 
 sup \int_{\omega}(|x-y| +| y|)d y \|f_0\|=  
$$
$$
 \frac{1}{2\tau_j\sigma_0^2} 
 sup (\int_{-b}^{0}(x-2y)d y + \int_{0}^x x d y +\int_{x}^{b}(2y-x)d y )
 \|f_0\|=  
 $$
 \begin{equation}
 \label{A1}
 \frac{1}{2\tau_j\sigma_0^2} 
 sup ( 2b^2 +x^2)  \|f_0\|= \frac{3 b^2}{2\tau_j \sigma_0^2} \|f_0\|  
 \end{equation}
 where $sup$ is over $x\in\omega = (-b,b)$ (or , equivalently, over 
  $x\in (0,b)$ ). 
  
  Similarly,
$$
\|A_{j2}f_1\|\leq 
$$
$$
 = (\frac{1}{2\sigma_0^2} \int_{\omega}| y|d y +
\sqrt{ \frac{\tau_j}{2}}\frac{1}{\sigma_0} 
sup ( e^{\frac{x^2}{2\tau_j\sigma_0^2}}
\int_{\omega}\int_{\frac{|x-y| +| y| }{\sigma _0\sqrt{2\tau_j}}}^{\infty} 
e^{-\theta ^2}d\theta)d y ) \|f_1\| \leq 
$$
 \begin{equation}
 \label{A2}
 ( \frac{ b^2}{2\sigma_0^2}+\sqrt{\frac{\pi \tau_j}{2}}\frac{b}{\sigma_0}) \|f_1\| . 
 \end{equation}
 
 To explain the last inequality we observe that
 $$
e^{X^2}\int_X^{\infty}e^{-\theta^2} d\theta = 
\int_0^{\infty}e^{-\sigma^2} \frac{\sigma}{\sqrt{X^2+\sigma^2}} d\sigma \leq
\int_0^{\infty}e^{-\sigma^2} d\sigma =\frac{\sqrt{\pi}}{2} 
 $$
 where we used the substitution $\theta=\sqrt{X^2+\sigma^2}$ and the
 well known probability integral. Using that $|x| \leq |x-y|+|y|$ and hence 
 $$
 \int_{\frac{|x-y| +| y| }{\sigma _0\sqrt{2\tau_j}}}^{\infty} 
e^{-\theta ^2}d\theta \leq
\int_{\frac{|x| }{\sigma _0\sqrt{2\tau_j}}}^{\infty} 
e^{-\theta ^2}d\theta
 $$
 and letting $X = \frac{x}{\sigma_0\sqrt{2\tau_j}}$ we yield \eqref{A2}.

From \eqref{systemd}, \eqref{A1}, \eqref{A2} by the triangle inequality we yield
$$
\|f_0\|-
\frac{1}{\tau_2-\tau_1}(\frac{\tau_2}{\tau_1}+\frac{\tau_1}{\tau_2})\frac{3b^2}{2\sigma_0^2}\|f_0\|-
$$
 $$
(\frac{\tau_1+\tau_2}{\tau_2-\tau_1}\frac{b^2}{2\sigma_0^2}+
\frac{\tau_2 \sqrt{\tau_1}+ \tau_1 \sqrt{\tau_2}}{\tau_2-\tau_1}
\frac{\sqrt{\pi}b}{\sqrt{2}\sigma_0})\|f_1\| \leq
\frac{1}{\tau_2-\tau_1}\|\tau_2 w_1-\tau_1 w_2\|
$$
and
$$
\|f_1\|-
\frac{1}{\tau_2-\tau_1}
(\frac{1}{\tau_1}+\frac{1}{\tau_2})\frac{3b^2}{2\sigma_0^2}\|f_0\|-
$$
 $$
\frac{1}{\tau_2-\tau_1}(\frac{b^2}{\sigma_0^2}+
\frac{\sqrt{\pi}b}{\sqrt{2}\sigma_0})(\sqrt{\tau_1}+\sqrt{\tau_2})\|f_1\| \leq
\frac{1}{\tau_2-\tau_1}\| w_1-w_2\|.
$$ 
Adding the first inequality and the second inequality multiplied by $\sqrt{\tau_1\tau_2}$ we obtain 
$$
\|f_0\|+ \sqrt{\tau_1\tau_2} \|f_1\|-
\frac{1}{(\tau_2-\tau_1)\tau_1\tau_2}
(\tau_1^2+\tau_2^2 +
\sqrt{\tau_1\tau_2}(\tau_1+\tau_2))\frac{3b^2}{2\sigma_0^2}\|f_0\|-
$$
$$
\frac{1}{\tau_2-\tau_1}
((\tau_1+\tau_2+2\sqrt{\tau_1\tau_2})\frac{b^2}{2\sigma_0^2}+
\sqrt{2\pi}(\sqrt{\tau_1}+\sqrt{\tau_2})\sqrt{\tau_1\tau_2}\frac{b}{\sigma_0})\|f_1\|\leq
$$
$$
\frac{1}{\tau_2-\tau_1}(\|\tau_2 w_1-\tau_1 w_2\|+\sqrt{\tau_1 \tau_2}\| w_1-w_2\|).
$$ 
 
Using the condition \eqref{uniqueness} we conclude that the factors
of $\|f_0\|, \|f_1\|$ on the left side are positive, so (by dividing by smallest of them
both sides) we arrive at \eqref{boundB}.

The proof is complete.

By Lemma 4.2 the linearized inverse option pricing problem is 
equivalent to the system of linear integral equations
$$
\int_{\omega}(K_0( x,y;\tau_1 ) f_0( y)+K_1( x,y;\tau_1 ) f_1( y)) dy = W_1(x),\quad
x\in \omega,  
$$
\begin{equation}
\int_{\omega}(K_0( x,y;\tau_2 ) f_0( y)+K_1( x,y;\tau_2 ) f_1( y)) dy = W_2(x),\quad
x\in \omega,  
\label{system1}
\end{equation}
where
$$
K_0( x,y;\tau) =S_*\int_{\frac{|x-y| +| y| }{\sigma _0\sqrt{2\tau}}}^{\infty} 
e^{-\tau ^2}d\tau,\; \; S_*=\frac{s^*}{\sigma^2_0\sqrt{\pi}}, 
$$
$$
K_1( x,y;\tau) =S_*(\frac{\sqrt{\tau}}{\sqrt{2}\sigma_0}(|y|-|x-y|)
e^{-\frac{(|x-y| +| y|)^2 }{2\tau \sigma_0^2}}+
$$
$$
(\frac{x^2-2xy}{\sigma_0^2}+\tau)
\int_{\frac{|x-y| +| y| }{\sigma _0\sqrt{2\tau}}}^{\infty} 
e^{-\tau ^2}d\tau). 
$$
and $W_j(x)=e^{-cx-d\tau_j}V_j(x)$ with $V_j(x)=U_j(x)-V_0(x,\tau_j)$
as defined in \eqref{final*}, \eqref{additionaldataV}.

Theorem 4.1 guarantees uniqueness and stability of a solution $(f_0, f_1)\in 
C(\bar\omega)\times C(\bar\omega)$ to this system under the condition 
(\ref{uniqueness}). It is not known whether this condition is necessary
 for uniqueness in the linearized inverse problem. It is not hard to
show that uniqueness for (\ref{system1}) holds under 
a similar condition for $(f_0,f_1)\in L^2(\omega) \times L^2(\omega)$.

At this point we do not have an explicit bound on $C$ in (\ref{boundB}).

We can not claim the existence of the solution of the linearized inverse option pricing problem, but we can show that under the conditions of Theorem 4.1 the range of the operator defined by the left side of \eqref{system1} (considered from $C(\bar\omega)\times C(\bar\omega)$ into $C^2(\bar\omega)\times C^2(\bar\omega)$) has the codimension 4. Indeed, for any $(W_1,W_2) \in C^2(\bar\omega)\times C^2(\bar\omega)$ the Fredholm system \eqref{system} (which follows from \eqref{system1}) has a unique solution $(f_0,f_1)\in C(\bar\omega)\times C(\bar\omega)$, due to uniqueness guaranteed by theorem 4.1. Obviously, using this $(f_0,f_1)$ in the left side of \eqref{system1}
produces the right side $C_1+C_2 x$ with some constant vectors $C_1, C_2$. 
Observe that we can not conclude that $C_1=0, C_2=0$, so while the system
(\ref{system}) follows from the system \eqref{system1}, it is not equivalent to this
equation. In any event, the range of the operator from left side of \eqref{system1} has the codimension not greater than 4   We will show that it is exactly 4.

 To do it we will observe that  the system (\ref{system1}) 
has no solution when $W_1, W_2$ are linear functions of $x$. Indeed, let $(f_0, f_1)$
be this solution. In the both cases $W_j^{\prime\prime}(x)=0$, so the solution $(f_0, f_1)$ to the equation (\ref{system1}) is zero due to uniqueness guaranteed by
Theorem 4.1.  But then $W_j=0$ and we have a contradiction.

At present we do not know an  exact description of the range of the (matrix)
linear integral operator defined by the left side of \eqref{system1}. 

For analytical and especially numerical purposes it can be useful to handle 
$f_j\in C(\mathbb{R})$. The results of section 4 can be extended to these functions which are constants on the intervals $(-\infty,-b], [b,+\infty)$. The analysis will be a little more technical and the condition \eqref{uniqueness} will be replaced by more  complicated one.

\section{Conclusion}

Uniqueness result of this paper can be certainly extended to the
many-dimensional case and probably to general parabolic equations with
variable coefficients. The inverse option pricing problem is a particular 
case of the more general inverse diffusion problem which has a 
probabilistic interpretation. We are not aware of any uniqueness results 
about recovery of space dependent diffusion rate from probability of distribution 
at a fixed moment of time, and moreover of a (special) time dependent rate
from probability at few moments of time. The method of this paper can be
 applied at least to a linearized version of this inverse probabilistic 
problem. Observe that while in case of one unknown space dependent coefficient
(scalar case) one can obtain complete uniqueness results \cite{Is} under certain monotonicity conditions on boundary data (which are not satisfied for the inverse option pricing problem) by using maximum or positivity principles, the (vectorial) case of several unknown spacial functions corresponds to systems of equations when maximum
principles are only in very special situations. There are almost no uniqueness results
in the vectorial case.

The proposed reconstruction algorithm is expected to perform very well
when volatility is not changing fast with respect to stock price
$s$ and is changing very slow with respect to time. Sudden
and dramatic changes of market situations most likely can not
be properly described by our model and more generally by
the Black-Scholes equation.  Probably, a minor 
modification of the proposed model (replacing ${\mathbb R}$ in (\ref{BS*y}) 
by a finite interval) can eliminate difficulties with existence theorem and generate 
even better numerical algorithms. Observe, that to
find continuous $f_j$ the data $U_j$ must be at least twice differentiable
on $\omega$, so the real market data are in need of a proper interpolation,
minimizing the size of second derivatives of $U_j$. A choice of an appropriate
smoothing interpolation and an intensive numerical testing will be
a subject of future work. Our immediate goal was
to develop some mathematical theory of the inverse option pricing problem  and
to suggest a numerical algorithm.
The next step is an additional testing on simulated and 
real data. We already made preliminary numerical tests which show very fast
convergence of iterative methods. The proposed method provides with simplest representation of the inverse problem and hence promises to be fastest and most stable of existing numerical methods.

For simplicity we considered only European options. We hope to adjust
the linearization technique to American and more complicated options,
which are in particular described by free boundary problems and multidimensional
parabolic equations. So far there are actually no analytic results in this very important practical case.

{\bf Aknowledgement:}   

   This research was in part supported by the NSF grant DMS 10-08902 and by Emylou Keith and Betty Dutcher Distinguished Professorship at the Wichita State University.

 \end{document}